\documentstyle[12pt]{article}

\def\Z2{Z\!\!\!\!\!\!Z_{\,2}}

\newcommand{\bn}{\begin{equation}}
\newcommand{\ed}{\end{equation}}
\newcommand{\bneqn}{\begin{eqnarray}}
\newcommand{\edeqn}{\end{eqnarray}}
\newcommand{\bnth}{\begin{theorem}}
\newcommand{\edth}{\end{theorem}}
\newcommand{\bnpr}{\begin{proposition}}
\newcommand{\edpr}{\end{proposition}}
\newcommand{\bnlm}{\begin{lemma}}
\newcommand{\edlm}{\end{lemma}}
\newcommand{\bndef}{\begin{definition}}
\newcommand{\eddef}{\end{definition}}
\newcommand{\bntb}{\begin{tabular}}
\newcommand{\edtb}{\end{tabular}}

\newtheorem{definition}{Definition}[section]
\newtheorem{lemma}{Lemma}[section]
\newtheorem{proposition}{Proposition}[section]
\newtheorem{theorem}{Theorem}[section]

\normalbaselines
\begin{document}
\title{ON A GENERAL ANALYTICAL FORMULA FOR $U_q(su(3))$-CLEBSCH-GORDAN
COEFFICIENTS}
\author{R.M. {ASHEROVA}\thanks{Institute of Physics and Power Engineering,
Obninsk; e-mail: asherova@nucl-th.sinp.msu.ru},
YU.F. SMIRNOV\thanks{e-mail: smirnov@nucl-th.sinp.msu.ru}
and V.N. TOLSTOY\thanks{e-mail: tolstoy@nucl-th.sinp.msu.ru}}
\date{}
\maketitle
\begin{center}
\vskip -10pt
{\large Institute of Nuclear Physics, Moscow State University\\[3pt]
119899 Moscow \& Russia 
}\\ [30pt]
\end{center}

\begin{abstract}
We present the projection operator method in combination with
the Wigner-Racah calculus of the subalgebra $U_q(su(2))$ for calculation
of Clebsch-Gordan coefficients (CGCs) of the quantum algebra $U_q(su(3))$.
The key formulas of the method are couplings of the tensor and projection
operators and also a tensor form for the projection operator of
$U_q(su(3))$. We obtain a very compact general analytical formula for
the $U_q(su(3))$ CGCs in terms of the $U_q(su(2))$ Wigner $3nj$-symbols.
\end{abstract}
\maketitle

\section{Introduction}
It is well known that the Clebsch-Gordan coefficients (CGCs) of the
unitary Lie algebra $u(n)$ ($su(n)$) have numerous applications in
various fields of theoretical and mathematical physics. For example,
many algebraic models of nuclear theory(interacting boson model (IBM),
Elliott $su(3)$ model, $su(4)$ supermultiplet scheme of Wigner,
the shell model, and so on) demand
the CGCs for $su(6)$, $su(5)$, $su(3)$, $su(4)$ and $su(n)$.
Analogously, in quark models of hadrons we need the CGCs of $su(3)$,
$su(4)$, etc. The theory of the $su(n)$ CGCs is connected with the
theory of special functions, combinatorial analysis, topology, etc.

There are several methods for the calculation of CGCs of $su(n)$
($u(n)$) and other Lie algebras:
recursion method; method of employment of explicit bases of
irreducible representations; method of generating invariants;
method of tensor operators, where the Wigner-Eckart theorem is used;
projection operator method; coherent state method; combined methods.

It is well known that the method of projection operators for usual
(non-quantized) Lie algebras \cite{AST1,T1} and superalgebras \cite{T1}
is powerful and universal method for a solution of many problems
in the representation theory. In particular, the method allows
to develop the detailed theory of Clebsch-Gordan
coefficients and another elements of Wigner-Racah calculus (including
compact analytic formulas of these elements and their symmetry
properties) \cite{PST1} and so on. It is evident that the projection
operators of quantum groups \cite{T2} play the same role in their
representation theory.

In this paper we present the projection operator method in combination
with the Wigner-Racah calculus of the subalgebra $U_q(su(2))$
\cite{STK1} for calculation of CGCs of the quantum algebra $U_q(su(3))$.
The key formulas of the method are couplings of the tensor and
projection operators and also a tensor form for the projection operator
of $U_q(su(3))$. It should be noted that the first application of
this method was for the $su(3)$ case in 
\cite{PST1}.
Some simple elements of this approach were also used in
\cite{AS} for the $U_q(su(n))$ case.
Also, the coherent state method in combination with the Wigner-Racah
calculus was applied in \cite{R} for $u(n)$.

\setcounter{equation}{0}
\section{Gelfand-Tsetlin basis} 

Let $\Pi:=\{\alpha_{1},\alpha_{2}\}$ be a system of simple roots of
the Lie algebra $sl(3)$ ($=\!sl(3,{\bf C})\simeq A_2$), endowed
with  the following scalar product:
$(\alpha_{1},\alpha_{1})\!=\!(\alpha_{2},\alpha_{2})\!=\!2$,
$(\alpha_{1},\alpha_{2})\!=\!(\alpha_{2},\alpha_{1})\!=\!-1$.
The root system $\Delta_+$ of $sl(3)$ consists of the roots
$\alpha_1,\alpha_1\!+\!\alpha_2,\alpha_2$.
The quantum Hopf algebra $U_{q}(sl(3))$ is generated by the Chevalley
elements $q^{\pm h_{\alpha_i}}$, $e_{\pm\alpha_i}$ $(i=1,2)$ with the
relations:
\begin{equation}
\begin{array}{rcl}
q^{h_{\alpha_i}}q^{-h_{\alpha_i}}&\!\!=\!\!&
q^{-h_{\alpha_i}}q^{h_{\alpha_j}}=1,\quad\;\;
q^{h_{\alpha_i}}q^{h_{\alpha_j}}=q^{h_{\alpha_j}}q^{h_{\alpha_i}},
\quad\;\;q^{h_{\alpha_i}}e_{\alpha_j}^{}q^{-h_{\alpha_i}}=
q^{(\alpha_i,\alpha_j)}e_{\alpha_j}^{},
\\[7pt]
[e_{\alpha_i}^{},e_{-\alpha_j}^{}]\!\!&=\!\!&
\delta_{ij}^{}\,[h_{\alpha_i}]~,
\qquad\quad
[[e_{\pm\alpha_i},e_{\pm\alpha_j}]_{q}^{},e_{\pm\alpha_j}]_{q}^{}=0
\qquad{\rm for}\;\;|i-j|=1~.
\label{GT1}
\end{array}
\end{equation}
Here and elsewhere we use the standard notation
$[a]:=(q^a-q^{-a})/(q-q^{-1})$, and
$[e_{\alpha}^{},e_{\beta}^{}]_q:=e_{\alpha}^{}e_{\beta}^{}-
q^{(\alpha,\beta)}e_{\beta}^{}e_{\alpha}^{}$.
The Hopf structure of $U_{q}(u(3))$ is given by
\begin{equation}
\begin{array}{rcccl}
\Delta_{q}(h_{\alpha_i})\!\!\!&=\!\!\!&
h_{\alpha_i}\otimes 1 +1\otimes h_{\alpha_i},
\qquad\qquad\qquad\quad
S_{q}(h_{\alpha_i})\!\!\!&=\!\!\!&-h_{\alpha_i},
\\[7pt]
\Delta_{q}(e_{\pm\alpha_i}^{})\!\!\!&=\!\!\!&e_{\pm\alpha_i}^{}\!
\otimes q^{\frac{1}{2}h_{\alpha_i}}
+q^{-\frac{1}{2}h_{\alpha_i}}\!\otimes e_{\pm\alpha_i}^{},
\qquad
S_{q}(e_{\pm\alpha_i}^{})\!\!\!&=\!\!\!&-q^{\pm1}e_{\pm\alpha_i}^{}.
\end{array}
\label{GT2}
\end{equation}
For construction of the composite root vectors
$e_{\pm(\alpha_1+\alpha_2)}^{}$ we fix the normal ordering in
$\Delta_+$: $\alpha_1,\;\alpha_1+\alpha_2,\;\alpha_2$.
According to this ordering
we put
\begin{equation}
e_{\alpha_1+\alpha_2}^{}:=[e_{\alpha_1}^{},e_{\alpha_2}^{}]_{q^{-1}}^{}~,
\qquad
e_{-\alpha_1-\alpha_2}^{}:=[e_{-\alpha_2}^{},e_{-\alpha_1}^{}]_{q}^{}~.
\label{GT3}
\end{equation}

Let us introduce another standard notations for the Cartan-Weyl
generators:
\begin{equation}
\begin{array}{rcccccl}
e_{12}^{}&\!\!:=\!\!&e_{\alpha_1}^{},\qquad\;\;
e_{21}^{}&\!\!:=\!\!&e_{-\alpha_1}^{},\qquad\;\;
e_{11}^{}-e_{22}^{}&\!\!:=\!\!&h_{\alpha_1}^{},
\\[5pt]
e_{23}^{}&\!\!:=\!\!&e_{\alpha_2}^{},\qquad\;\;
e_{32}^{}&\!\!:=\!\!&e_{-\alpha_2}^{},\qquad\;\;
e_{22}^{}-e_{33}^{}&\!\!:=\!\!&h_{\alpha_2}^{},
\\[5pt]
e_{13}^{}&\!\!:=\!\!&e_{\alpha_1+\alpha_2}^{},\quad \;
e_{31}^{}&\!\!:=\!\!&e_{-\alpha_1-\alpha_2}^{},\quad\;
e_{11}^{}-e_{33}^{}&\!\!:=\!\!&h_{\alpha_1}^{}\!+h_{\alpha_2}.
\end{array}
\label{GT4}
\end{equation}
The explicit formula for the extremal projector for the quantum groups
\cite{T2} specialized to the case of $U_q(sl(3))$ has the form
\begin{equation}
p\bigl(U_q(sl(3))\bigr)=p_{12}^{}p_{13}^{}p_{23}^{}~,
\label{GT5}
\end{equation}
where the elements $p_{ij}$ ($1\le i<j\le3$) are given by
\begin{equation}
p_{ij}^{}=\sum\limits_{n=0}^{\infty}
\frac{(-1)^n}{[n]!}\,\varphi_{ij,n}^{}\,e_{ij}^{n}e_{ji}^{n}~,\quad\;\;
\varphi_{ij,n}^{}=q^{-(j-i-1)n}\Bigr\{\prod\limits_{s=1}^{n}
[e_{ii}^{}\!-e_{jj}^{}\!+j-i+s]\Bigr\}^{-1}\!.
\label{GT6}
\end{equation}
The extremal projector $p:=p\bigl(U_q(sl(3))\bigr)$ satisfies
the relations:
\begin{equation}
e_{ij}^{}p=pe_{ji}^{}=0 \quad (i<j)~,
\qquad p^{2}=p~.
\label{GT7}
\end{equation}
The quantum algebra $U_{q}(su(3))$ can be considered as the
quantum algebra $U_{q}(sl(3))$ endowed with the additional
Cartan involution $(^*)$:
\begin{equation}
h_{\alpha_i}^*=h_{\alpha_i}~,\qquad
e_{\pm\alpha_i}^*=e_{\mp\alpha_i}^{}~,\qquad
q^*=q\;({\rm or}\;\,\bar{q}:=q^{-1})~.
\label{GT8}
\end{equation}

Let $(\lambda\mu)$ be a finite-dimensional irreducible representation
(IR) of $U_q(su(3))$ with the highest weight $(\lambda\mu)$
($\lambda$ and $\mu$ are nonnegative integers).
The vector of the highest weight, denoted by the symbol
$\bigr|(\lambda\mu)h\bigl>$, satisfy the relations
\begin{equation}
h_{\alpha_1}\bigr|(\lambda\mu)h\bigl> 
=\lambda\bigr|(\lambda\mu)h\bigl>~,\quad\;
h_{\alpha_2}\bigr|(\lambda\mu)h\bigl>=
\mu\bigr|(\lambda\mu)h\bigl>~,\quad\;
e_{ij}^{}\bigr|(\lambda\mu)h\bigl> 
=0\quad (i<j)~.
\label{GT9}
\end{equation}
Labelling of another basis vectors in IR $(\lambda\mu)$ depends upon
choice of subalgebras of $U_q(su(3))$ (or in another words, depends
upon which reduction chain from $U_q(su(3))$ to subalgebras is chosen).
Here we use the Gelfand-Tsetlin reduction chain:
\begin{equation}
U_q(su(3))\supset U_q(u_{Y}^{}(1))\otimes
U_q(su_{T}^{}(2))\supset U_q(u_{T_{0}}^{}(1))~,
\label{GT10}
\end{equation}
where the subalgebra $U_q(su_{T}^{}(2))$ is generated by the elements
\begin{equation}
T_{+}:=e_{23}^{}~,\qquad T_{-}:=e_{32}^{}~,\qquad
T_{0}:=\mbox{\large$\frac{1}{2}$}(e_{22}^{}-e_{33}^{})~,
\label{GT11}
\end{equation}
the subalgebra $U_q(u_{T_0}^{}(1))$ is generated by $q^{T_{0}}$,
and 
$U_q(u_{Y}^{}(1))$ is generated by $q^{Y}$ 
(In the classical (non-deformed) case in the elementary
particle theory the subalgebra $su_{T}^{}(2)$ is called the T-spin
algebra and the element ${Y}$ is the hypercharge operator), where:
\begin{equation}
Y=-\mbox{\large$\frac{1}{3}$}\bigl(2h_{\alpha_1}^{}+
h_{\alpha_2}^{}\bigr)~.
\label{GT12}
\end{equation}
In the case of the reduction chain (\ref{GT10}) the basis vectors
of IR $(\lambda\mu)$ are denoted by
\begin{equation}
\bigl|(\lambda\mu)jtt_{z}\bigr>~.
\label{GT13}
\end{equation}
Here the 
set $jtt_{z}$ characterize the hypercharge $Y$ and the T-spin and
its projection:
\begin{equation}
\begin{array}{rcl}
q^{T_0}\bigl|(\lambda\mu)jtt_{z}\bigr>&\!\!=\!\!&q^{t_z}
\bigl|(\lambda\mu)jtt_{z}\bigr>~,\qquad
q^{Y}\bigl|(\lambda\mu)jtt_{z}\bigr>=q^{y}
\bigl|(\lambda\mu)jtt_{z}\bigr>~,
\\[9pt]
T_{\pm}\bigl|(\lambda\mu)jtt_{z}\bigr>&\!\!=\!\!&\sqrt{[t\mp t_z]
[t\pm t_z\!+\!1]}\bigl|(\lambda\mu)jtt_{z}\!\pm\!1\bigr>~,
\end{array}
\label{GT14}
\end{equation}
where the parameter $j$ is connected with the eigenvalue $y$ of
the operator $Y$ as follows
$y=-\mbox{\large$\frac{1}{3}$}\bigl(2\lambda+\mu\bigr)+2j$.
It is not hard to show that the orthonormalized vectors (\ref{GT13})
can be represented in the following form
\begin{equation}
\bigr|(\lambda\mu)jtt_{z}\bigr>=N^{(\lambda\mu)}_{\,jt}
P^{\,t}_{\!t_{z};t}\;e_{31}^{j+\frac{1}{2}\mu-t}
e_{21}^{j-\frac{1}{2}\mu+t}\bigr|(\lambda\mu) h\bigr>~,
\label{GT15}
\end{equation}
where $P^{\,t}_{\!t_{z};t_z'}$ is the general projection operator
of the quantum algebra $U_q(su_{T}^{}(2))$ \cite{STK1}, and
the normalizing factor $N^{(\lambda\mu)}_{jt}$ has the form
\begin{equation}
N^{(\lambda\mu)}_{\,jt}=\left(\frac{[\lambda+\frac{1}{2}\mu-j+t+1]!
[\lambda+\frac{1}{2}\mu-j-t]![\frac{1}{2}\mu+j+t+1]!
[\frac{1}{2}\mu-j+t]!}{q^{2j+\mu-2t}\,[\lambda]![\mu]!
[\lambda+\mu+1]![j+\frac{1}{2}\mu-t]!
[j-\frac{1}{2}\mu+t]![2t+1]!}\right)^{\frac{1}{2}}\!.
\label{GT16}
\end{equation}
The quantum numbers $jt$ are taken all nonnegative integers and
half-integers such that the sum $\mbox{\large$\frac{1}{2}$}\mu+j+t$
is an integer and they are subjected to the constraints:
\begin{equation}
\left\{\begin{array}{rcl}
\frac{1}{2}\mu+j-t&\!\!\ge\!\!&0~,\quad\;-\frac{1}{2}\mu+j+t\ge0~,
\\[3pt]
\frac{1}{2}\mu-j+t&\!\!\ge\!\!&0~,
\qquad\frac{1}{2}\mu+j+t\ge\lambda+\mu~.
\end{array}
\right.
\label{GT17}
\end{equation}
For every fixed $t$ the projection $t_z$ runs the values
$t_z=-t,-t+1,\ldots,t-1,t$. These results can be obtained
from the explicit form of the Gelfand-Tsetlin bases for the case
$U_q(su(n))$ \cite{T2} specializing to the given case $U_q(su(3))$.

\setcounter{equation}{0}
\section{Couplings of tensor and projection operators}
Let $\{R_{j_z}^{j(q)}\}$ be an irreducible tensor operator (ITO)
of the rank $j$, that is $(2j+1)$-components $R_{j_z}^{j(q)}$ are
transformed with respect to the $U_q(su_T(2))$ adjoint action as
the $U_q(su_T(2))$ basis vectors $\bigr|jj_z\bigr>$ of the spin $j$:
\begin{equation}
T_{i}^{}\triangleright R_{j_z}^{j(q)}:=({\rm ad}_q\,T_{i}^{})
R_{j_z}^{j(q)}\equiv\bigl(({\rm id}\otimes S_q)
\Delta_q(T_{i}^{})\bigr)\circ R_{j_z}^{j(q)}=\sum\limits_{j'_z}
\bigl<jj'_z\bigl|T_{i}^{}\bigr|jj_z\bigr>R_{j'_z}^{j(q)}~,
\label{CTP1}
\end{equation}
where $(a\otimes b)\circ x=axb$. The tensor operator of the type
$\{R_{j_z}^{j(q)}\}$ will be also called the left irreducible tensor
operators (LITO) because the generators $T_{i}$ ($i=\pm,0$) act to
the left-side of the components $R_{j_z}^{j(q)}$.
(The given denotation of the ITOs is differed from one
of the papers \cite{STK1} by the replacement $q$ by $q^{-1}$).
Following to the paper \cite{PST1} we also introduce a right
irreducible tensor operator (RITO) denoted by the tilde symbol
$\{\tilde{R}_{j_z}^{j(q)}\}$, on which the $U_q(su_T(2))$
generators $T_i$ act on the right-side, namely
\begin{equation}
T_{i}^{}\triangleleft\tilde{R}_{j_z}^{j(q)}:=
({\rm ad}_q^*\,T_{i}^{})\,\tilde{R}_{j_z}^{j(q)}\equiv
\tilde{R}_{j_z}^{j(q)}\stackrel{\leftarrow}{\circ}
\bigl((\tilde{S}_q\otimes{\rm id})\tilde{\Delta}_q(T_{i}^*)\bigr)=
\sum\limits_{j'_z}\bigl<jj'_z\bigl|T_{i}^{}\bigr|jj_z\bigr>
\tilde{R}_{j'_z}^{j(q)}~.
\label{CTP2}
\end{equation}
where $x\stackrel{\leftarrow}{\circ}(a\otimes b)=axb$, and
$\tilde{\Delta}_q$ is the opposite coproduct
($\tilde{\Delta}_q=\Delta_{\bar{q}}$) and $\tilde{S}_q$ is
the corresponding antipode ($\tilde{S}_q=S_{\bar{q}}$).
It is not hard to verify that any LITO $\{R_{j_z}^{j(q)}\}$
is the RITO $\{\tilde{R}_{j_z}^{j(q)}\}$:
$R_{j_z}^{j(q)}=(-1)^{j_z}q^{j_z}\tilde{R}_{-j_z}^{j(q)}$.
The projection operator set $\{P^{t}_{\!t_z;t'_z}\}$ for a fixed IR
$t$ and for various $t_z$ and $t'_z$  will be called the
${\bf P}^{t}$-operator. It is not hard to see that the subset of
the left components of this operator satisfy the relations for
the LITO ${\bf R}^{j(q)}\!:=\!\{R_{j_z}^{j(q)}\}$ if we understand
the action "$\triangleright$" of the generator $T_i$ as the usual
multiplication of the operators $T_i$ and $P^{t}_{\!t_z;t'_z}$ and
the subset of the right components of the ${\bf P}^{t}$-operator
satisfy the relations for the RITO
$\tilde{\bf R}^j\!:=\!\{\tilde{R}_{j_z}^{j(q)}\}$ if we understand the
action $"\triangleleft"$ as the usual multiplication of the operators
$P^{t}_{\!t_z;t'_z}$ and $T_i^{*}$: 
\begin{eqnarray}
T_i^{}\triangleright P^{t}_{\!t_z;t'_z}:=T_i^{}P^{t}_{\!t_z:t'_z}
=\sum\limits_{t'\!'_{\!z}}\bigl<tt''_z\bigl|T_i^{}\bigr|tt_z\bigr>
P^{t}_{\!t'\!'_{\!z};t'_z}~,
\label{CTP3}
\\[0pt]
T_i^{}\triangleleft P^{t}_{\!t_z;t'_z}:=P^{t}_{\!t_z:t'_z}T_i^*=
\sum\limits_{t'\!'_{\!z}}\bigl<tt''_z\bigl|T_i^{}\bigr|tt'_z\bigr>
P^{t}_{\!t_z;t'\!'_{\!z}}~.
\label{CTP4}
\end{eqnarray}
Using the $U_q(su_T(2))$ CGCs we can couple the LITO ${\bf R}^{j(q)}$
with the left components of the ${\bf P}^{t}$-operator and the RITO
$\tilde{\bf R}^{j(q)}$ with the right components of
the ${\bf P}^{t}$-operator:
\begin{eqnarray}
{t'\atop t'_z}\!\!\Big\{{\bf R}^{j(q)}\mathop{\stackrel{.}{\otimes}}
{\bf P}^{t}_{\!;t_z}\!\Big\}_{\!q}:=
\sum\limits_{j_zt'\!'_{\!z}}
\bigl(j j_z\,tt''_z\bigl|t't'_z\bigr)_{q}\,R_{j_z}^{j(q)}
P^{t}_{\!t'\!'_{\!z},t_z}~,
\label{CTP5}
\\[0pt]
\Big\{{\bf P}^{t}_{t_z;}\mathop{\stackrel{.}{\otimes}}
{\bf\tilde{R}}^{j(q)}\Big\}_{\!q}{\!t'\atop t'_z}:=
\sum\limits_{j_zt'\!'_{\!z}}\bigl(j j_z\,tt''_z\bigl|t't'_z\bigr)_{q}\,
P^{t}_{\!t_z,t'\!'_{\!z}}\tilde{R}_{j_z}^{j(q)}~.
\label{CTP6}
\end{eqnarray}
Here the symbol $\mathop{\stackrel{.}{\otimes}}$ means that we first
take the usual tensor product and then in a resulting expression we
replace the tensor product by the usual operator product.
It is not hard to show that the couplings (\ref{CTP5}) and (\ref{CTP6})
are connected as
follows
\begin{equation}
{\bf I\!R}^{j(q)}_{tt_z;t't'_z}:=\sqrt{[2t+1]}\;
{t\atop t_z}\!\!\Big\{{\bf R}^{j(q)}\mathop{\stackrel{.}{\otimes}}
{\bf P}^{t'}_{\!;t'_z}\Big\}_{\!q}=(-1)^{t'-t}\,
\sqrt{[2t'+1]}\;\Big\{{\bf P}^{t}_{t_z;}\mathop{\stackrel{.}{\otimes}}
{\bf\tilde{R}}^{j(q)}\!\Big\}_{\!q}{\!\!t'\atop t'_z}~.
\label{CTP7}
\end{equation}
Using (\ref{CTP7}) and an unitary relation of the $U_q(su(2))$ CGCs
\cite{STK1} one can obtain the following useful permutation relations
between the components of the tensors ${\bf R}^{j(q)}$,
$\tilde{\bf R}^{j(q)}$ and ${\bf P}^{t}$-operator:
\begin{eqnarray}
R_{j_z}^{j(q)}P^{t}_{\!t_z;t'_z}=\sum\limits_{t'\!'t'\!'_{\!z}}\,
(-1)^{t-t''}\,\sqrt{\frac{[2t+1]}{[2t''+1]}}\;
\bigl(jj_z\,tt_z\bigl|t''t''_z\bigr)_{q}\;
\Big\{{\bf P}^{t'\!'}_{\!t'\!'_{\!z};}
\mathop{\stackrel{.}{\otimes}}
\tilde{\bf R}^{j(q)}\!\Big\}_{\!q}{\!\!t\atop t'_z}~,
\label{CTP8}
\\[5pt]
P^{t}_{\!t_z;t'_z}\tilde{R}_{j_z}^{j(q)}=
\sum\limits_{t'\!'t'\!'_{\!z}}\,
(-1)^{t-t''}\,\sqrt{\frac{[2t+1]}{[2t''+1]}}\;
\bigl(jj_z\,tt'_z\bigl|t''t''_z\bigr)_{q}\;
{t\atop t_z}\!\!\Big\{{\bf R}^{j(q)}\!
\mathop{\stackrel{.}{\otimes}}
{\bf P}^{t'\!'}_{\!;t'\!'_{\!z}}\Big\}_{\!q}~.
\phantom{a}
\label{CTP9}
\end{eqnarray}
We can show that the monomials $e_{21}^{n}e_{31}^{m}$ and
$e_{12}^{n}{e'}_{\!\!13}^{m}$ are components of ITOs with respect
to the adjoint action of the subalgebra $U_q(su_T(2))$:
\begin{eqnarray}
R^{j(q)}_{j_z}=
\sqrt{\frac{[2j]!}{[j-j_{z}]![j+j_{z}]!}}\;q^{2j^2-j}\,
e_{21}^{j+j_{z}}e_{31}^{j-j_{z}}q^{-jh_{\alpha_1}-(j-j_z)T_0},
\label{CTP10}
\\[7pt]
{R'}^{j(q)}_{\!j_z}=
\sqrt{\frac{[2j]!}{[j-j_{z}]![j+j_{z}]!}}\;q^{-2j^2+j}\,
e_{12}^{j-j_{z}}{e'}_{\!\!13}^{j+j_{z}}q^{-jh_{\alpha_1}-(j+j_z)T_0}.
\label{CTP11}
\end{eqnarray}
where the generator ${e'}_{\!\!13}$ is defined according to the
inverse normal ordering: $\alpha_2,\;\alpha_1+\alpha_2,\;\alpha_1$,
i.e. ${e'}_{\!\!13}=[e_{23},e_{12}]_{q^{-1}}$.
These ITOs have the remarkable properties:
{\it A result of the coupling of two ITOs of the type (\ref{CTP10})
or (\ref{CTP11}) is non-zero only for an irreducible component of
the maximal rank, e.g.}
\begin{equation}
\Bigl\{{\bf R}^{j(q)}\!\mathop{\stackrel{.}{\otimes}}
{\bf R}^{j'(q)}\!\bigr\}_{\!q}\!{\!\!j''\atop j''_z}=
\delta_{j'\!'\!,j+j'}\,
R^{j+j'(q)}_{\,j'\!'_{\!z}}.
\label{CTP12}
\end{equation}
The property is also useful in applications:
{\it For ITOs of the type (\ref{CTP12}) the relation is valid}
\begin{equation}
R^{j(q)}_{j_z}{\bf I\!R}^{j'(q)}_{tt_z;t't'_z}=
\sum\limits_{t'\!'t'\!'_{\!z}}
\sqrt{\frac{[2t+1]}{[2t''+1]}}
\;\bigl(jj_z\,tt_z\bigl|t''t''_z\bigr)_{q}\,
U\!\bigl(jj't''t';j\!+\!j't\bigr)_{\!q}\,
{\bf I\!R}^{j+j'(q)}_{t'\!'t'\!'_{\!z};t't'_z}.
\label{CTP13}
\end{equation}
Here $U\bigl(\ldots;\ldots\bigr)_{\!q}$ is recoupling coefficient
which can be expressed via the stretched $q$-$6j$-symbols
of $U_q(su_T(2))$ \cite{STK1}:
\begin{equation}
U\!\bigl(jj't''t';j\!+\!j't\bigr)_{\!q}=
(-1)^{j+j'+t'+t'\!'}\sqrt{[2j+2j'+1][2t+1]}\;
{\displaystyle\left\{{j\atop t'}\;
{j'\atop t''}\;{j\!+\!j'\atop t}\right\}_{\!q}}.
\label{CTP14}
\end{equation}
Using (\ref{CTP9}) and (\ref{CTP10}) we can present the basis
vectors (\ref{GT13}) in the form of
\begin{equation}
\bigl|(\lambda\mu)jtt_z\bigr>={\cal F}_{\!-}^{(\lambda\mu)}(jtt_z)
\bigl|(\lambda\mu)h\bigr>={\cal N}^{(\lambda\mu)}_{jt}\,
{\bf I\!R}^{j(q)}_{tt_z;\frac{1}{2}\mu\frac{1}{2}\mu}
\bigl|(\lambda\mu)h\bigr>~.
\label{CTP15}
\end{equation}
The normalizing factor ${\cal N}^{(\lambda\mu)}_{jt}$ is given by
\begin{equation}
\begin{array}{rcl}
{\cal N}^{(\lambda\mu)}_{jt}&\!\!=\!\!& 
(-1)^{2j}q^{(j+\frac{1}{2}\mu-t)(j-\frac{1}{2}\mu+t)+j\lambda+
\frac{1}{2}\mu(j+\frac{1}{2}\mu-t)-2j^2+j+t-\frac{1}{2}\mu}
\\[5pt]
&&{}\times
\mbox{\large$\sqrt{\frac{[j-\frac{1}{2}\mu+t]![j+\frac{1}{2}\mu-t]!}
{
[2j]![\mu+1]}}$}\,
\bigl(j\,\mbox{\large$\frac{1}{2}$}\mu\!-\!t\;tt\bigl|
\mbox{\large$\frac{1}{2}$}\mu\mbox{\large$\frac{1}{2}$}\mu\bigr)_{q}
\;N^{(\lambda\mu)}_{jt}~.
\end{array}
\label{CTP16}
\end{equation}
With the help of (\ref{CTP15}) we easy find the action of the ITO
(\ref{CTP10}) on the Gelfand-Tsetlin basis:
\begin{equation}
R^{j'(q)}_{j'_z}\bigl|(\lambda\mu)jtt_z\bigr>=
\sum\limits_{t'\!'t'\!'_{\!z}}
\;\bigl(j'j'_z\,tt_z\bigl|t''t''_{\!z}\bigr)_{q}\,
\bigl<(\lambda\mu)j''t''\bigr\|R^{j'(q)}\bigl\|
(\lambda\mu)jt\bigr>_{\!q}\,\bigl|(\lambda\mu)j''t''t''_z\bigr>~,
\label{CTP17}
\end{equation}
where
\begin{equation}
\bigl<(\lambda\mu)j''t''\bigr\|R^{j'(q)}\bigl\|(\lambda\mu)jt\bigr>_{\!q}
=\delta_{j''\!,\,j+j'}^{}\sqrt{\frac{[2t+1]}{[2t''+1]}}\;
\frac{{\cal N}^{(\lambda\mu)}_{jt}}{{\cal N}^{(\lambda\mu)}_{j''t''}}\,
U\bigl(j'j\,t''\,\mbox{\large$\frac{1}{2}$}\mu;j'\!\!+\!jt\bigr)_{\!q}~.
\label{CTP18}
\end{equation}

\setcounter{equation}{0}
\section{Tensor form of the 
projection operator}
It is obvious that the extremal projector (\ref{GT5}) can be
presented in the form
\begin{equation}
p(U_q(su(3))=p(U_q(su^{}_T(2))
\bigl(p_{12}^{}p_{13}^{}\bigr)p(U_q(su^{}_T(2)).
\label{TFP1}
\end{equation}
Now we present the middle part of
(\ref{TFP1}) in the terms of the $U_q(su_T(2))$ tensor operators
(\ref{CTP10}) and (\ref{CTP11}). To this end, we substitute the
explicit expression (\ref{GT6}) for the factors $p_{12}$ and
$p_{13}$, and combine monomials  $e_{21}^{n}e_{31}^{m}$ and
$e_{12}^{n}e_{13}^{m}$.
After some summation manipulations we obtain the following expression
for the extremal projection operator $p\!:=\!p(U_q(su(3)))$ in terms
of the tensor operators (\ref{CTP10}) and (\ref{CTP11}):
\begin{equation}
p=p(U_q(su^{}_T(2))
\Bigr(\sum\limits_{jj_z}
A_{jj_z}
\tilde{R}^{j(q)}_{j_z}{R'}^{j(q)}_{\!j_z}\Bigl)p(U_q(su^{}_T(2)).
\label{TFP2}
\end{equation}
Here
\begin{equation}
A_{jj_z}=\frac{(-1)^{3j}[\varphi_{12}^{}]
[\varphi_{12}^{}+j+j_z-1]![\varphi_{13}^{}]!}
{[2j]![\varphi_{12}^{}+2j]![\varphi_{13}^{}+j+j_z]!}\,
q^{4j^2+j+2jh_{\alpha_1}+2(j+j_z)T_0},
\label{TFP3}
\end{equation}
where $\varphi_{1i+1}:=e_{11}\!-e_{i+1i+1}\!+i$ ($i=1,2$).
Below we assume that the $U_q(su(3))$ extremal projection operator $p$
acts in a weight space with the weight $(\lambda\mu)$ and in this
case the symbol $p$ is supplied with the index $(\lambda\mu)$,
$p^{(\lambda\mu)}$, and all the Cartan elements $h_{\alpha_i}$ on the
right side of (\ref{TFP2}) are  replaced by the corresponding weight
components $\lambda$ and $\mu$. Now we multiple the projector
$p^{(\lambda\mu)}$ from the left side by the lowering operator
${\cal F}_{\!-}^{(\lambda\mu)}(jtt_z)$ and from the right side by the
rising operator $\bigl({\cal F}_{\!-}^{(\lambda\mu)}(jtt_z)\bigr)^*$,
and by applying a relation of type (\ref{CTP13}) we finally find
the tensor form of the general $U_q(su(3))$ projection operator:
\begin{equation}
P^{(\lambda\mu)}_{\!\!jtt_z;j't't'_z}\!=
\sum\limits_{j''t''}B_{j''t''}^{(\lambda\mu)}\;
{\bf I\!R}^{j+j''(q)}_{tt_z,t''t''}\;
{\bf I\!R'}^{j''\!+j'(q)}_{\!t''t'',t't'_z}~,
\label{TFP4}
\end{equation}
were the coefficients $B_{j'\!'t'\!'}^{(\lambda\mu)}$ are given by
\begin{equation}
\begin{array}{l}
B_{j''t''}^{(\lambda\mu)}=
\mbox{\large$\frac{(-1)^{2j+j'\!+j''\!-t'+t''}q^{\phi}
[\lambda+1][\mu+1][\lambda+\mu+2]}
{[\lambda+\frac{1}{2}\mu+j''\!+t''\!+2]!
[\lambda+\frac{1}{2}\mu+j''\!-t''+1]![2j'']!}$}\;
{\displaystyle\left\{{j\atop t''}\;
{j''\atop t}\;{j\!+\!j''\atop\frac{1}{2}\mu}\right\}_{\!q}}\;
{\displaystyle\left\{{j'\atop t''}\;
{j''\atop t'}\;{j'\!+\!j''\atop\frac{1}{2}\mu}\right\}_{\!q}}
\\[17pt]
{}\;\;\times\mbox{\large$\left(\frac{
[\lambda+\frac{1}{2}\mu-j+t+1]!
[\lambda+\frac{1}{2}\mu-j-t]!
[\lambda+\frac{1}{2}\mu-j'\!+t'\!+1]!
[\lambda+\frac{1}{2}\mu-j'\!-t']![2j+2j''\!+1][2j'+2j''\!+1]}
{[2j]![2j']![2t+1][2t'+1]}
\right)^{\frac{1}{2}}$}\!,
\end{array}
\label{TFP5}
\end{equation}
\begin{equation}
{}\quad\phi=\varphi(\lambda,\mu,j,t)+\varphi(\lambda,\mu,j'\!,t')
-2\varphi(\lambda,\mu,j''\!,t'')+j''(4\lambda\!+\!2\mu\!-\!1)
+4t''\!-\!2\mu\!-\!3j'.
\label{TFP6}
\end{equation}
Here and elsewhere we use the notation $\varphi(\lambda,\mu,j,t):=
\frac{1}{2}(\frac{1}{2}\mu+j-t)(\frac{1}{2}\mu+j+t-3)+j(\lambda-2j+\!1)$.

\setcounter{equation}{0}
\section{General form of 
Clebsch-Gordan coefficients}
For convenience we introduce the short notations:
$\Lambda:=(\lambda\mu)$ and $\gamma:=jtt_{z}$ and therefore the
basis vector $\bigl|(\lambda\mu)jtt_z\bigr>$ will be is denoted by
$\bigl|\Lambda\gamma\bigr>$.
Let $\{|\Lambda_i\gamma_i\bigr>\}$ be bases of two IRs
$\Lambda_i$ $(i=1,2)$. Then
$\{|\Lambda_1\gamma_1\bigr>|\Lambda_2\gamma_2\bigr>\}$ be
a basis in the representation
$\Lambda_1\otimes \Lambda_2$ of $U_q(su(3))\otimes U_q(su(3))$.
In this representation there is an another coupled basis
$|\Lambda_1\Lambda_2\!:s\Lambda_3\gamma_3\bigr>_q$ with respect to
$\Delta_q(U_q(su(3)))$ where the index $s$ classifies multiple
representations $\Lambda_3$.We can expand the coupled basis in terms
of the tensor uncoupled basis
$\{|\Lambda_1\gamma_1\bigr>|\Lambda_2\gamma_2\bigr>\}$:
\begin{equation}
\bigl|\Lambda_1\Lambda_2\!:s\Lambda_3\gamma_3\bigr>_q=
\sum_{\gamma_1,\gamma_2}^{}\bigl(\Lambda_1\gamma_1\,
\Lambda_2\gamma_2\bigl|s\Lambda_3\gamma_3\bigr)_{\!q}\,
\bigr|\Lambda_1\gamma_1\bigr>\bigr|\Lambda_2\gamma_2\bigr>~,
\label{CGC1}
\end{equation}
where the matrix element
$\bigl(\Lambda_1\gamma_1\,\Lambda_2\gamma_2|s\Lambda_3\gamma_3\bigr)_{\!q}$
is the Clebsch-Gordan coefficient of $U_q(su(3))$.
In just the same way as for the non-quantized Lie algebra $su(3)$
(see \cite{PST1}) we can show that any CGC of $U_q(su(3))$
can be represented in terms of the linear combination of the matrix
elements of the projection operator (\ref{TFP4})
\begin{equation}
\bigl(\Lambda_1\gamma_1\,\Lambda_2\gamma_2|s\Lambda_3\gamma_3\bigr)_q=
\sum_{\gamma_2'}^{}C(\gamma_2')\,\bigl<\Lambda_1\gamma_1\bigr|
\bigl<\Lambda_2\gamma_2\bigr|\Delta_q(P_{\!\gamma_3,h}^{\Lambda_3})
\bigr|\Lambda_1h\bigl>\bigr|\Lambda_2\gamma_2'\bigr>~.
\label{CGC2}
\end{equation}
Classification of multiple representations $\Lambda_3$ in the
representation $\Lambda_1\otimes\Lambda_2$ is special problem and we
shall not touch it here. For the non-deformed algebra $su(3)$ this
problem was considered in details in \cite{PST1}. Concerning of
the matrix elements in the right-side of (\ref{CGC2}) we give here
an explicit expression for the more general matrix element:
\begin{equation}
\bigl<\Lambda_1\gamma_1\bigr|\bigl<\Lambda_2\gamma_2\bigr|
\Delta_q(P_{\!\gamma_3,\gamma_3'}^{\Lambda_3})
\bigr|\Lambda_1\gamma_1'\bigr>\bigr|\Lambda_2\gamma_2'\bigr>~.
\label{CGC3}
\end{equation}
Using (\ref{TFP4}) and the Wigner-Racah calculus for the subalgebra
$U_q(su(2))$ \cite{STK1} (analogously to the non-quantized Lie algebra
$su(3)$ \cite{PST1}) it is not hard to obtain the following result:
\begin{equation}
\begin{array}{lcl}
\mbox{\normalsize$\bigl<\Lambda_1\gamma_1\bigr|
\bigl<\Lambda_2\gamma_2\bigr|
\Delta_q(P_{\gamma_3,\gamma_3'}^{\Lambda_3})
\bigr|\Lambda_1\gamma_1'\bigr>\bigr|\Lambda_2\gamma_2'\bigr>$}=
\mbox{\normalsize$\bigl(t_1t_{1z}\,t_2t_{2z}\bigr|t_3t_{3z}\bigr)_{q}\;
\bigl(t_1t_{1z}'\,t_2t_{2z}'\bigr|t_3't_{3z}'\bigr)_{q}$}
\\[14pt]
\qquad\;\;{}\times
[\lambda_3+1][\mu_3+1][\lambda_3+\mu_3+2]\;\;
A\!\sum\limits_{j_1'\!'j_2'\!'t_1'\!'t_2'\!'t_3'\!'}
C_{j_1'\!'j_2'\!'t_1'\!'t_2'\!'t_3'\!'}
\\[17pt]
\qquad\;\;{}\times
\left\{\!\!\!\!\begin{array}{cccc}
&j_1\!\!-\!j_1'\!'&j_2\!\!-\!j_2'\!'&
j_1\!\!+\!j_2\!\!-\!j_1'\!'\!\!-\!j_2'\!'\\
&t_1'\!'&t_2'\!'&t_3'\!'\\
&t_1&t_2&t_3\\
\end{array}
\!\right\}_{\!q}
\left\{\!\!\!\!\begin{array}{cccc}
&j_1'\!\!-\!j_1'\!'&j_2'\!\!-\!j_2'\!'&
j_1'\!\!+\!j_2'\!\!-\!j_1'\!'\!\!-\!j_2'\!'\\
&t_1'\!'&t_2'\!'&t_3'\!'\\
&t_1'&t_2'&t_3'
\end{array}
\!\right\}_{\!q}.
\end{array}
\label{CGC4}
\end{equation}
Here
\begin{equation}
\begin{array}{rcl}
A&\!\!=\!\!& 
\mbox{\large$\left(\frac{[2t_1+1][2t_2+1][2j_1+1]![2j_2+1]!
[\lambda_3+\frac{1}{2}\mu_3-j_3+t_3+1]!
[\lambda_3+\frac{1}{2}\mu_3-j_3-t_3]!}
{[\lambda_1+\frac{1}{2}\mu_1-j_1+t_1+1]!
[\lambda_1+\frac{1}{2}\mu_1-j_1-t_1]!
[\lambda_2+\frac{1}{2}\mu_2-j_2+t_2+1]!
[\lambda_2+\frac{1}{2}\mu_2-j_2-t_2]![2j_3]!}
\right.$}
\\[17pt]
&&\times\mbox{\large$\left.
\frac{[2t_1'+1][2t_2'+1]
[2j_1'+1]![2j_2'+1]!
[\lambda_3+\frac{1}{2}\mu_3-j_3'+t_3'+1]!
[\lambda_3+\frac{1}{2}\mu_3-j_3'-t_3']!}
{[\lambda_1+\frac{1}{2}\mu_1-j_1'+t_1'+1]!
[\lambda_1+\frac{1}{2}\mu_1-j_1'-t_1']!
[\lambda_2+\frac{1}{2}\mu_2-j_2'+t_2'+1]!
[\lambda_2+\frac{1}{2}\mu_2-j_2'-t_2']![2j'_3]!}
\right)^{\frac{1}{2}}$},
\end{array}
\label{CGC5}
\end{equation}
\begin{equation}
\begin{array}{lcl}
C_{j_1'\!'j_2'\!'t_1'\!'t_2'\!'t_3'\!'}=
\mbox{\large$\frac{
(-1)^{2(j_1\!+j_2\!+j_3'\!-j_1'\!'\!-j_2'\!')}\,q^{\psi}
[2(j_1+j_2-j''_1-j''_2)+1]![2(j_1'+j_2'-j''_1-j''_2)+1]!
}{[2j''_1]![2j''_2]![2j_1-2j''_1]![2j_2-2j''_2]!
[2j_1'-2j''_2]![2j_2'-2j''_2]!
[2(j_1\!+\!j_2-\!j_3-\!j''_1\!-\!j''_2)]!}$}
\\[17pt]
{}\;\times\mbox{\large$\frac{
[\lambda_1\!+\frac{1}{2}\mu_1\!-\!j''_1\!+t''_1\!+\!1]!
[\lambda_1\!+\frac{1}{2}\mu_1\!-\!j''_1\!-t''_1]!
[\lambda_2\!+\frac{1}{2}\mu_2\!-\!j''_2\!+t''_2\!+\!1]!
[\lambda_2\!+\frac{1}{2}\mu_2\!-\!j''_2\!-t''_2]!
[2t''_1\!+\!1][2t''_2\!+\!1][2t''_3\!+\!1]}{
[\lambda_3+\frac{1}{2}\mu_3+j_1\!+\!j_2-\!j_3\!-\!
j''_1\!-\!j''_2\!+t''_3\!+\!2]!
[\lambda_3+\frac{1}{2}\mu_3+j_1\!+\!j_2-\!j_3\!-\!
j''_1\!-\!j''_2\!-\!t''_3\!+\!1]!}$}
\\[17pt]
\qquad{}\times
{\displaystyle\left\{{j_1\!\!-\!\!j''_1\atop\frac{1}{2}\mu_1}\;
{j''_1\atop t_1}\;{j_1\atop t''_1}\right\}_{\!q}}\;
{\displaystyle\left\{{j_2\!\!-\!\!j''_2\atop\frac{1}{2}\mu_2}\;
{j''_2\atop t_2}\;{j_2\atop t''_2}\right\}_{\!q}}\;
{\displaystyle\left\{{j_3\atop t''_3}\;
{\,{j_1\!\!+\!\!j_2\!\!-\!\!j_3\!\!-\!\!j''_1\!\!-\!\!j''_2\atop t_3}}\;
{\,{j_1\!\!+\!\!j_2\!\!-\!\!j''_1\!\!-\!\!j''_2\atop\frac{1}{2}\mu_3}}
\right\}_{\!q}}
\\[21pt]
\qquad{}\times
{\displaystyle\left\{{j_1'\!\!-\!\!j''_1\atop\frac{1}{2}\mu_1}\;
{j''_1\atop t_1'}\;{j_1'\atop t''_1}\right\}_{\!q}}\;
{\displaystyle\left\{{j_2'\!\!-\!\!j''_2\atop\frac{1}{2}\mu_2}\;
{j''_2\atop t_2'}\;{j_2'\atop t''_2}\right\}_{\!q}}\;
{\displaystyle\left\{{j_3'\atop t''_3}\;{\,{j_1'\!\!+
\!j_2'\!\!-\!\!j_3'\!\!-\!\!j''_1\!\!-\!\!j''_2\atop t_3'}}\;
{\,{j_1'\!\!+\!\!j_2'\!\!-\!\!j''_1\!\!-\!\!j''_2
\atop\frac{1}{2}\mu_3}}\right\}_{\!q}},
\end{array}
\label{CGC6}
\end{equation}
where $\psi=\sum_{i=1}^{2}\Bigl(2\varphi(\lambda_i,\mu_i,j''_i,t''_i)-
\varphi(\lambda_i,\mu_i,j_i,t_i)-\varphi(\lambda_i,\mu_i,j'_i,t'_i)
-t_i(t_i+1)-t_i'(t'_i+1)\Big) \\[2pt]
-2\varphi(\lambda_3,\mu_3,j_3''\!,t_3'')
+\varphi(\lambda_3,\mu_3,j_3,t_3)
+\varphi(\lambda_3,\mu_3,j'_3\!,t_3')
+j_3''(4\lambda_3\!+2\mu_3+2)-2t_3''(t''_3\!-1)-2\mu_3\\[2pt]
-(j_2\!+\!j_2'\!-\!2j''_2)(2\lambda_1\!+\!\mu_1\!-6j_1'')
+4(j_1\!-\!j_1'')(j_2\!-\!j_2'')+4(j'_1\!-\!j_1'')(j'_2\!-\!j_2'')
-(j_3\!+\!j''_3)(j_3\!+\!j''_3\!+\!1)\\[2pt]
-(j'_3\!+\!j''_3)(j'_3\!+\!j''_3\!+\!1)$,
$\;j''_3:=j_1\!+\!j_2\!-\!j_3\!-\!j_1''\!-\!j_2''=
j'_1\!+\!j'_2\!-\!j'_3\!-\!j_1''\!-\!j_2''$.\\[-23pt]

\section*{Acknowledgments} 
This work was supported by
RFBR-99-01-01163, and  by the program of
French-Russian scientific cooperation (CNRS grant PICS-608 and
grant RFBR-98-01-22033, V.N. Tolstoy).\\[-23pt]





\begin{thebibliography}{99}

\bibitem{AST1}
R.M. Asherova, Yu.F. Smirnov, and V.N. Tolstoy,
Teor. Mat. Fiz. {\bf 8}, no. 2, {255} (1971);
Teor. Mat. Fiz. {\bf 15}, no. 1, {107} (1973);
Mat. Zamet.  {\bf 26}, {15} (1979).

\bibitem{T1}
V.N. Tolstoy,
Uspekhi Mat. Nauk {\bf 40}, no. 4 (244), {225} (1985);
Uspekhi Mat. Nauk {\bf 44}, no. 1 (265), {211} (1989),
[transl. in Russian Math. Surveys. {\bf 44}, no. 1, {257} (1989)].

\bibitem{PST1}
Z. Pluhar, Yu.F. Smirnov, and V.N. Tolstoy, (1981)
Charles University preprint, 
(Prague, 1981); J. Phys. A: Math. Gen. {\bf 19}, no. 1, {21} (1986).

\bibitem{T2}
V.N. Tolstoy,
Lectures Notes in Phys. {\bf 370},  {118} (Springer, Berlin, 1990).

\bibitem{STK1}
Yu.F. Smirnov, V.N. Tolstoy, and Yu.I. Kharitonov,
Yad. Fiz. {\bf 53}, no. 4, {959} (1991),
[transl. in Soviet J. Nucl. Phys. {\bf 53}, no. 4, {593} (1991)],
Yad. Fiz. {\bf 53}, no. 6, {1746} (1991),
[transl. in Soviet J. Nucl. Phys. {\bf 53}, no. 6, {1068} (1991)];
Yad. Fiz. {\bf 54}, no. 3, {721} (1991),
[transl. in Soviet J. Nucl. Phys. {\bf 54}, no. 3, {437} (1991)];
Yad. Fiz. {\bf 55}, no. 10, {2863} (1992),
[transl. in Soviet J. Nucl. Phys. {\bf 55}, 
{1599} (1992)];
Yad. Fiz. {\bf 56}, no. 5, {223} (1993).

\bibitem{AS}
S. Alisauskas, and Yu.F. Smirnov, J. Phys. A {\bf 27}, {5925} (1994).

\bibitem{R}
D.J. Rowe, and J. Repka, J. Math. Phys. {\bf 37}, {6530} (1997).

\end{thebibliography}
\end{document}